\documentclass[reqno]{amsart}
\usepackage{hyperref}
\usepackage{amsmath}
\usepackage{amssymb}
\usepackage{amscd}
\usepackage{graphics}
\usepackage{latexsym}
\usepackage{stmaryrd}

\theoremstyle{plain}
\newtheorem{theorem}{Theorem}[section]
\newtheorem{thm}[theorem]{Theorem}

\newtheorem{lem}[theorem]{Lemma}

\theoremstyle{definition}
\newtheorem{defn}[theorem]{Definition}

\newtheorem{rmk}[theorem]{Remark}

\newtheorem{notat}[theorem]{Notation}

\theoremstyle{remark}

\newcommand{\QQ}{\mathbb{Q}}

\newcommand{\ZZ}{\mathbb{Z}}

\newcommand{\PP}{\mathbb{P}}

\newcommand{\mc}{\mathcal}
\newcommand{\mf}{\mathfrak}

\newcommand{\OO}{\mc{O}}

\newcommand{\SP}{\text{Spec }}

\newsavebox{\sembox}
\newlength{\semwidth}
\newlength{\boxwidth}

\newsavebox{\semrbox}
\newlength{\semrwidth}
\newlength{\boxrwidth}

\title{Degenerations of rationally connected varieties and PAC fields} 
\author[Starr]{Jason Starr} 
\begin{document}


\begin{abstract}
A degeneration of a separably rationally connected variety over a
field $k$ contains a geometrically irreducible subscheme if $k$
contains the algebraic closure of its subfield.  If $k$ is a perfect
PAC field, the degeneration has a $k$-point.  This generalizes
~\cite[Theorem 21.3.6(a)]{FriedJarden}: a degeneration of a Fano
complete intersection over $k$ has a $k$-point if $k$ is a perfect PAC
field containing the algebraic closure of its prime subfield.
\end{abstract}


\maketitle


\section{Statement of results} \label{sec-sor}

Recently, a number of fields long known to be $C_1$ were proved to
satisfy an \emph{a priori} stronger property related to rationally
connected varieties.
\begin{enumerate}
\item[(i)] 
Every rationally connected variety defined over the
function field of a curve over a characteristic 0 algebraically closed
field has a 
rational point, ~\cite{GHS}.
\item[(ii)] 
Every separably rationally connected
variety defined over the function field of a curve over an
algebraically closed field of arbitrary characteristic has a closed
point, ~\cite{dJS}. 
\item[(iii)]
Every smooth, rationally chain connected
variety over a finite field has a rational point, ~\cite{Esnault}.  
\end{enumerate}
Moreover, in
each of these cases, degenerations of these varieties also have
rational points, at least under some mild hypotheses on the
degeneration.  

This article considers the same problem for perfect PAC fields
containing an algebraically closed field.  Such fields are known to be
$C_1$, ~\cite[Theorem 21.3.6(a)]{FriedJarden}.  The main theorem is
the following.

\begin{thm} \label{thm-main}
Let $k$ be a perfect PAC field containing the algebraic closure of its
prime subfield.  Let $X_k$ be a $k$-scheme which is the closed fiber
of a proper, flat scheme over a DVR whose geometric generic fiber is
separably rationally connected (in the sense of ~\cite{dJS}).
Then $X_k$ has a $k$-point.  
\end{thm}

\begin{rmk} \label{rmk-main}
This gives a new proof of ~\cite[Theorem 21.3.6(a)]{FriedJarden},
i.e., every perfect PAC field containing an algebraically closed field
is $C_1$.   
\end{rmk}

This should be compared to the following theorems of Koll\'ar
and de Jong respectively.

\begin{thm} ~\cite{KAx} \label{thm-K}
Let $k$ be a characteristic $0$ PAC field.  Let $X_k$ be a $k$-scheme
whose base-change $X\otimes_k \overline{k}$ is the closed fiber of a
projective, flat scheme over a DVR whose geometric generic fiber is a
Fano manifold. 
Then
$X_k$ has a $k$-point.  In particular, $k$ is $C_1$.  
\end{thm}

\begin{thm}[de Jong] \label{thm-dJ}
Let $k$ be a characteristic $0$ field having a point in every
rationally connected $k$-scheme and containing $\overline{\QQ}$.  
Let $X_k$ be a $k$-scheme which is the closed fiber
of a proper, flat scheme over a DVR whose geometric generic fiber is
rationally connected.  
Then $X_k$ has a $k$-point.
In particular, $k$ is $C_1$. 
\end{thm}

Theorem ~\ref{thm-main} is a consequence of the following
more precise result.

\begin{thm} \label{thm-main2}
Let $k$ be a field containing the algebraic closure of its prime
subfield.  Let $X$ be a proper $k$-scheme which is the closed fiber of
a proper, flat scheme over a DVR whose geometric generic fiber is 
separably rationally connected (in the sense of ~\cite{dJS}).
There exists a
closed subscheme $Y$ of $X$ such that $Y\otimes_k \overline{k}$ is
irreducible.
\end{thm}



\medskip\noindent
\textbf{Acknowledgments.}  I thank A. J. de Jong for explaining his
proof of Theorem ~\ref{thm-dJ} some years ago.  I thank J\'anos
Koll\'ar for encouraging me to write-up this article and pointing out
~\cite[Theorem 21.3.6(a)]{FriedJarden}.




\section{Proofs} \label{sec-pfs}


\begin{defn} \label{defn-int}
Let $R$ be a DVR and let $X_R$ be an $R$-scheme.
A \emph{finite type model} of $(R,X_R)$
consists of a datum
$$
((P,D)\rightarrow (S,\mf{s}),\SP R \rightarrow P, X_P)
$$
of 
\begin{enumerate}
\item[(i)]
a Dedekind domain $S$ of finite type
over $\ZZ$ or over $\QQ$
and a maximal ideal $\mf{s}$ of $S$, 
\item[(ii)]
a normal, flat, projective
$S$-scheme $P$ with geometrically integral generic fiber
and an integral Weil divisor $D$ of $P$ contained in $P_\mf{s}$,
\item[(iii)]
a dominant morphism $\SP R \rightarrow P$, 
\item[(iv)]
and a $P$-scheme $X_P$
\end{enumerate} 
such
that the inverse image of $D$ equals $\mf{m}_R$ and the base change of
$X_P$ equals $X_R$.
\end{defn}

\begin{lem} \label{lem-int}
If $X_R$ is a finite type $R$-scheme, there exists a finite type model 
such that $X_P$ is of finite type over $P$.  

If $X_R$
is proper, resp. projective, there exists a finite type model such
that $X_P$ is proper, resp. projective.  
\end{lem}

\begin{proof}
This follows easily from ~\cite[\S 8]{EGA3}, results about
normalization, e.g., ~\cite[Scholie 7.8.3]{EGA4}, and Nagata
compactification, ~\cite{Ncomp}.
\end{proof}

Let $B$ be a scheme. 
As in ~\cite{Jou}, for each integer $e$ denote by $\text{Gr}(e,N)$ the
Grassmannian over $B$ parametrizing codimension-$e$ linear subspaces
of fibers of $\PP^N_B$.  
Denote by $\Lambda_e$
the universal codimension-$e$ linear subscheme of
$\text{Gr}(e,N)\times_B \PP^N_B$.  

\begin{notat} \label{notat-GZ}
For every $B$-scheme $T$ and $B$-morphism $i:T\rightarrow \PP^N_B$,
denote by $\mc{Z}_{i,e}$ the fiber product $T\times_{\PP^N_B}
\Lambda_e$.  If $i$ is inclusion of a subscheme, denote this by
$\mc{Z}_{T,e}$.  Observe that for a $k$-morphism $j:U\rightarrow T$,
$\mc{Z}_{i\circ j,e}$ equals $U\times_T \mc{Z}_{i,e}$.  
Denote by $\text{pr}_G:\mc{Z}_{i,e}\rightarrow \text{Gr}(e,N)$ and
$\text{pr}_T: \mc{Z}_{i,e} \rightarrow T$ the 2 projections.  

Denote
by $G_{i,e}$ the maximal open subscheme of $\text{Gr}(e,N)$ such that
$G_{i,e}\times_{\text{Gr}(e,N)} \mc{Z}_{i,e} \rightarrow G_{i,e}$ is
flat.  Denote by $\mc{Z}_{i,e,G}\rightarrow G_{i,e}$ the base change
of $\text{pr}_G$ to $G_{i,e}$.
\end{notat}  

Let $k$ be a field and 
let $I$ be a geometrically irreducible, 
locally closed subscheme of
$\PP^N_k$.  Denote by $d+1$ the dimension of $I$.
Denote $G_{I,d}$ by $G$ and denote $\mc{Z}_{I,d,G}$ by $\mc{Z}$.
Denote by $\text{pr}_G:\mc{Z}\rightarrow G$ and
$\text{pr}_I:\mc{Z}\rightarrow I$ the 2 projections.   

Let $H$ be a $k$-scheme and let $f:H\rightarrow G$ be a
morphism.  Denote by $\mc{Z}_H$ the fiber product $H\times_G \mc{Z}$.
Denote by $\text{pr}_H:\mc{Z}_H \rightarrow H$ and
$\text{pr}_\mc{Z}:\mc{Z}_H \rightarrow \mc{Z}$ the 2 projections.
Denote by $p:\mc{Z}_H \rightarrow I$ the composition $\text{pr}_I\circ
\text{pr}_\mc{Z}$.  

\begin{lem} \label{lem-irr}
If the dimension of $I$ is positive, if
$H$ is geometrically irreducible and if $f$ is dominant, then $\mc{Z}_H$ is
geometrically irreducible, and the geometric generic fiber of $p$ 
is irreducible. 
\end{lem}

\begin{proof}
By definition of $G$, $\mc{Z}$ is flat over $G$ of relative dimension
$1$.  
By ~\cite[Th\'eor\`eme I.6.10.2]{Jou}, the geometric generic fiber of
$\mc{Z}\rightarrow G$ is irreducible.  Since $f$ is dominant,
$\mc{Z}_H \rightarrow H$ is flat of relative dimension $1$ with
irreducible geometric generic fiber.  Since $H$ is also geometrically
irreducible, $\mc{Z}_H$ is geometrically irreducible.

Because $f$ is dominant and $\text{pr}_I$ are dominant, 
also $\text{pr}_\mc{Z}$, and thus the composition $p$, 
is dominant. 
Thus, to prove the geometric generic
fiber of $p$ is irreducible, 
it suffices to prove
the geometric generic fiber of 
$$
\text{pr}_1:\mc{Z}_H \times_I \mc{Z}_H \rightarrow \mc{Z}_H
$$ 
is irreducible.  

The morphism $\text{pr}_H\circ \text{pr}_1:\mc{Z}_H \times_I \mc{Z}_H
\rightarrow H$ is the base-change of $\mc{Z}_{p,d} \rightarrow
\text{Gr}(d,N)$ by $H\rightarrow \text{Gr}(d,N)$.  The image of
$p$ is dense in $I$, thus has dimension $d+1$.  Thus, 
by
~\cite[Th\'eor\`eme I.6.10.2]{Jou}, the geometric generic fiber of
$\text{pr}_G:\mc{Z}_{p,d} \rightarrow \text{Gr}(d,n)$ is
irreducible.  Since $H\rightarrow \text{Gr}(d,N)$ is dominant, also
the geometric generic fiber of $\text{pr}_H\circ \text{pr}_1$ is
irreducible.  Denote the geometric generic fiber by $X$.  

Denote by
$Y$ the geometric generic fiber of $\text{pr}_H:\mc{Z}_H \rightarrow
H$.  By the argument above, $Y$ is also irreducible.  The morphism
$\text{pr}_1$ induces a dominant morphism $q:X\rightarrow Y$.  The
geometric generic fiber of $p$ is irreducible if and only if the
geometric generic fiber of $q$ is irreducible.  

Denote by $\widetilde{q}:\widetilde{Y} \rightarrow Y_\text{red}$ the
integral closure of $Y_\text{red}$ in
the fraction field of $X_\text{red}$.  By the Noether normalization
theorem, $\widetilde{q}$ is finite.  And the geometric generic fiber
of $X_\text{red}\rightarrow \widetilde{Y}$ is integral.
Thus the 
geometric generic fiber of $p$ is irreducible if and only if the
geometric generic fiber of $\widetilde{q}$ is
irreducible.

Because $X$ is irreducible, also $\widetilde{Y}$ is irreducible.  The
diagonal morphism
$\Delta:\mc{Z}_H \rightarrow \mc{Z}_H \times_I \mc{Z}_H$ gives a
section of $\text{pr}_1$.  This induces a section $s:Y\rightarrow X$
of $q$.  This induces a section $\widetilde{s}:Y_\text{red} \rightarrow
\widetilde{Y}$ of $\widetilde{q}$.  

Since $\widetilde{q}$ is finite, $\text{dim}(\widetilde{Y})$ equals
$\text{dim}(Y_\text{red})$.  
The image of $\widetilde{s}$ is an
irreducible closed subscheme of $\widetilde{Y}$ whose dimension equals
$\text{dim}(Y_\text{red})$, i.e., $\text{dim}(\widetilde{Y})$.
Therefore, the
image is an irreducible component of $\widetilde{Y}$.  Since
$\widetilde{Y}$ is irreducible, the image of $\widetilde{s}$ is all of
$\widetilde{Y}$.  Thus $\widetilde{s}$ is an inverse of
$\widetilde{q}$.  Therefore the geometric generic fiber of
$\widetilde{q}$ is a point, which is irreducible.
\end{proof}

\begin{proof}[Proof of Theorem~\ref{thm-main2}]
Let $R$ be a DVR whose residue field $k$ contains the algebraic
closure of its prime subfield.  Let $X_R$ be a proper $R$-scheme whose
geometric generic fiber is separably rationally connected.
By Lemma ~\ref{lem-int}, there exists a finite type model with
$X_P$ proper.  By hypothesis, $k$ contains the algebraic closure of
the residue field $\overline{\kappa(\mf{s})}$.  

If $P$ has relative dimension $0$ over $S$, then $P$ equals $S$.  Thus
$X_k$ is the base change of the proper $\kappa(\mf{s})$-scheme
$X_\mf{s}$.  Since $\kappa(\mf{s})$ is a finite extension of the prime
subfield of $k$, the algebraic closure
$\overline{\kappa(\mf{s})}$ is contained in $k$.  The base change
$Y$ of any
$\overline{\kappa(\mf{s})}$-point of $X_\mf{s}$ is a geometrically
irreducible subvariety of $X_k$.  

Thus, assume the relative dimension of $P$ over $S$ is positive, $d+1$.
Let $P\hookrightarrow \PP^N_S$ be a closed immersion.  Denote by
$G_{P,d}$ and $\mc{Z}_{P,d,G}$ the schemes over $S$ from
Notation ~\ref{notat-GZ}.  By definition, the geometric generic fiber
of $P$ over $S$ is normal.  In particular it is smooth in codimension
$1$ points.  Therefore, by ~\cite[Th\'eor\`eme I.6.10.2]{Jou}, the
geometric generic fiber of $\mc{Z}_{P,d,G}\rightarrow G_{P,d}$ is
smooth.   

By ~\cite[Theorem 6.1]{Artin}, there exists
an algebraic space $\Pi$ separated and locally of finite type over
$G_{P,d}$, and a universal morphism
$$
\sigma: \Pi\times_{G_{P,d}} \mc{Z}_{P,d,G} \rightarrow X_P
$$
whose composition with projection to $P$ equals the composition
$\text{pr}_P\circ \text{pr}_{\mc{Z}_{P,d,G}}$.  Since the geometric
generic fiber of $X_R$ is separably rationally connected, also the
geometric generic fiber of $X_P \rightarrow P$ is separably rationally
connected. 
Therefore, by ~\cite{dJS}, the geometric generic fiber of $\Pi
\rightarrow G_{P,d}$ is nonempty.  Therefore some irreducible
component $\Pi_i$ of $\Pi$ dominates $G_{P,d}$.

By ~\cite[Corollary 5.20]{Kn},
there exists a dense open subspace $U$ of $\Pi_i$ which is an 
affine scheme. 
There
exists a dense open immersion of $U$ in a projective $G_{P,d}$-scheme
$\overline{U}$.  The morphism $\sigma$ induces a rational
transformation 
$$
\sigma_U:\overline{U} \times_{G_{P,d}} \mc{Z}_{P,d,G} \dashrightarrow X_P.
$$
Denote by $\overline{V}$ the normalization of
$\overline{U}\times_{G_{P,d}} \mc{Z}_{P,d,G}$.  Denote by $V$ the
maximal open subscheme of $\overline{V}$ over which $\sigma_U$ extends
to a regular morphism.  By the valuative criterion of properness, the
complement of $V$ has codimension $2$.  In particular, some
irreducible component of $V_{\mf{s}}$ dominates $(G_{P,d})_\mf{s}$.

Since $k$ contains $\overline{\kappa(\mf{s})}$ and the function field
$\kappa(D)$ of $D$, 
it contains the function field $\kappa(I)$ of one of the irreducible components
$I$ of $D\otimes_{\kappa(\mf{s})}
\overline{\kappa(\mf{s})}$.
Let $H$ denote an irreducible component of
$V_{\mf{s}}\otimes_{\kappa(\mf{s})} \overline{\kappa(\mf{s})}$
dominating $(G_{P,d})_\mf{s}\otimes_{\kappa(\mf{s})}
\overline{\kappa(\mf{s})}$.  
Let $G$, $\mc{Z}$, etc., be as in Lemma~\ref{lem-irr}.
Denote by $\mc{Z}_k$ the base change $\mc{Z} \otimes_{\OO_I} k$.  Note
that the base change $X_I \otimes_{\OO_I} k$ equals $X_k$ by
definition of a finite type model.  

The transformation $\sigma_U$ determines
a morphism 
$$
\sigma_H: \mc{Z}_H \rightarrow X_I.
$$
The base change to $k$ is a morphism of $k$-schemes,
$$
\sigma_k: \mc{Z}_k \rightarrow X_k.
$$
Denote the closure of the image by $Y$.
By Lemma ~\ref{lem-irr}, $\mc{Z}_H$ is
geometrically irreducible and the geometric generic fiber of $\mc{Z}_H
\rightarrow I$ is irreducible.  Thus $\mc{Z}_k \otimes_k
\overline{k}$ is irreducible.  Therefore $Y$ is a closed subscheme of
$X_k$ such that $Y\otimes_k \overline{k}$ is irreducible.
\end{proof}

\begin{proof}[Proof of Theorem ~\ref{thm-main}]
If $k$ is a perfect PAC field, then every geometrically irreducible
$k$-scheme has a $k$-point.  In particular, $Y$ has a $k$-point.
Therefore $X_k$ has a $k$-point.
\end{proof}

Every field $k$ is the closed fiber of a DVR $R$ whose generic fiber has
characteristic $0$.  Every complete intersection in $\PP^n_k$ is
the closed fiber of a complete intersection in $\PP^n_R$ whose generic
fiber is smooth.  If the complete intersection satisfies the $C_1$
inequality, the generic fiber is a Fano manifold.  By ~\cite{KMM92c},
~\cite{Ca}, a Fano manifold in characteristic $0$ is rationally
connected.  Therefore Theorem ~\ref{thm-main} implies the complete
intersection in $\PP^n_k$ has a $k$-point if $k$ is a perfect PAC
field containing the algebraic closure of its prime subfield.  In
other words, every perfect PAC field containing an algebraically
closed field is $C_1$, cf. ~\cite[Theorem 21.3.6(a)]{FriedJarden}.

\bibliography{my}
\bibliographystyle{alpha}

\end{document}